\documentclass[12pt,noamsfonts]{amsart}

\usepackage{amsfonts}
\usepackage{amssymb}
\usepackage{amscd}

\newcommand{\marginlabel}[1]%
  {\mbox{}\marginpar{\raggedleft\hspace{0pt}\bfseries\sf#1}}

\def\NN{{\mathbb N}}

\def\CC{{\mathbb C}}
\def\AA{{\mathbb A}}

\def\PP{{\mathbf P}}

\newtheorem{lemma}{Lemma}[section]
\newtheorem{theorem}[lemma]{Theorem}

\newtheorem{proposition}[lemma]{Proposition}

\theoremstyle{definition}

\theoremstyle{remark}
\newtheorem{remark}[lemma]{Remark}

\begin{document}

\title{The log canonical threshold of 
homogeneous affine hypersurfaces}

\author[L. Ein]{Lawrence~Ein}
\address{Department of Mathematics, University of Illinois at Chicago,
851 South Morgan St., M/C. 249, Chicago, IL 60607-7045}
\email{{\tt ein@math.uic.edu}}

\author[M. Musta\c{t}\v{a}]{Mircea~Musta\c{t}\v{a}}
\address{Department of Mathematics, University of California,
Berkeley, CA, 94720 and Institute of Mathematics of the Romanian
Academy}
\email{{\tt mustata@math.berkeley.edu}}

\maketitle
\markboth{L.~Ein and M.~Musta\c{t}\v{a}}{The log canonical threshold}

\section{Introduction}

In \cite{CP}, I.~Cheltsov and J.~Park studied the log canonical threshold of
 singular hyperplane sections of complex smooth,
projective hypersurfaces.
Let $X\subset\PP^n$, $n\geq 4$, be a complex 
smooth hypersurface of degree $d$ and $Y$
a hyperplane section of $X$ (which has to be irreducible and reduced).
I.~Cheltsov and J.~Park
proved that $Y$ has isolated singular points and
they studied the log canonical threshold of the pair $(X,Y)$.
They showed that 
$$c(X,Y)\geq\min\{(n-1)/d, 1\},$$
and they conjectured that if $d=n$, then equality holds if and only if
$Y$ is a cone over a (smooth) hypersurface in some $\PP^{n-2}$. 
Moreover, they showed that their conjecture follows from the Minimal 
Model Program.

The purpose of this note is to generalize their
lower bound in the case of an arbitrary hypersurface in $\PP^{n-1}$
and to give a direct proof of their conjecture in our more general setting.
In fact, we get these results also for the affine cone over the hypersurface.
 The main ingredient is  
the description of the log canonical threshold in terms of the
asymptotic growth of the jet schemes from \cite{Mu2}.

\medskip

Here are our results. Let $Y\subset\PP^{n-1}$, $n\geq 2$
be a complex hypersurface of degree $d\geq 1$ 
and let $Z\subset\AA^n$ be the affine cone over $Y$. Suppose that
$\dim\,{\rm Sing}\,(Z)=r$.

\begin{theorem}\label{main1}
With the above notation, we have the following lower bound
for the log canonical threshold of $(\AA^n,Z)${\rm :}
$$c(\AA^n,Z)\geq\min\{(n-r)/d, 1\}.$$
\end{theorem}

\begin{theorem}\label{main2}
If in addition $d\geq n-r+1$, then $c(\AA^n,Z)=(n-r)/d$ if and only if
$Z=T\times\AA^r$, for some hypersurface $T$ 
 in an $(n-r)$-dimensional affine space.
\end{theorem}

\begin{remark}
Since we have 
$$c(\PP^{n-1}, Y)=c(\AA^n\setminus\{0\}, Z\setminus\{0\})\geq c(\AA^n, Z),$$
it follows from Theorem~\ref{main1} that 
$c(\PP^{n-1}, Y)\geq\min\{(n-r)/d, 1\}$. Moreover, if
$d\geq n-r+1$ and $c(\PP^{n-1},Y)=(n-r)/d$, we deduce from Theorem~\ref{main2}
that $Y$ is the projective cone (with a $\PP^{r-1}$ vertex)
over a (smooth) hypersurface in some
$\PP^{n-r-1}$. 

The converse is well-known: if $Y\subset\PP^{n-1}$
 is the projective cone 
 over a smooth
hypersurface 
of degree $d\geq n-r+1$ in $\PP^{n-r-1}$, then $c(\PP^{n-1}, Y)=(n-r)/d$.
However,
for completeness, we will include an argument for this assertion
in the spirit of this paper in Proposition~\ref{standard} below.
\end{remark}

In order to make the connection between the way we stated our results
and the results in \cite{CP}, we make the following

\begin{remark}
Suppose we are in the situation in \cite{CP}: $X\subset\PP^n$ is a smooth
hypersurface and $Y=X\cap H$, where $H\subset\PP^n$ is a hyperplane.
We have the following equality:
$$c(H,Y)=c(X,Y).$$
See, for example, Theorem~\ref{ingred} below for justification.

Since Cheltsov and Park proved that in this case $r\leq 1$, 
Theorems~\ref{main1} and \ref{main2} give in particular their lower bound
and their conjectured characterization of equality.
\end{remark}

\section{Jet scheme dimension computations}

For the standard definition of the log canonical threshold,
as well as for equivalent definitions in singularity theory, we refer to
\cite{Ko}. We will take as definition the characterization from
\cite{Mu2} in terms of jet schemes.

Recall that for an arbitrary scheme $W$
(of finite type over $\CC$), the $m$th jet scheme $W_m$
is a scheme of finite type over $\CC$ characterized by
$${\rm Hom}({\rm Spec}\,A, W_m)\simeq{\rm Hom}({\rm Spec}\,A[t]/(t^{m+1}),
W),$$
for every $\CC$-algebra $A$. Note that $W_m(\CC)=
{\rm Hom}({\rm Spec}\,\CC[t]/(t^{m+1}), W),$ and in fact, we will be
interested only in the dimensions of these spaces.

\begin{theorem}\label{ingred}{\rm (\cite{Mu2} 3.4)}
If $X$ is a smooth, connected variety of dimension $n$, and $D\subset X$
is an effective divisor, then the log canonical threshold of $(X,D)$
is given by
$$c(X,D)=n-\sup_{m\in\NN}\frac{\dim\,D_m}{m+1}.$$
Moreover, there is $p\in\NN$ such that $c(X,D)=n-(\dim\,D_m)/(m+1)$
whenever $p\mid (m+1$.
\end{theorem}

\medskip

It is easy to write down equations for jet schemes.
We are interested in the jet schemes of a hypersurface
$Z\subset\AA^n$ defined by a polynomial $F\in\CC[X_i;1\leq i\leq n]$.
The jet scheme $Z_m$ is a subscheme of $\AA^{(m+1)n}={\rm Spec}\,R_m$,
where $R_m=
\CC[X_i,X'_i,\ldots,X_i^{(m)};i]$.
If $D:R_m\longrightarrow R_{m+1}$ is the unique $\CC$-derivation 
such that $D(X_i^{(j)})=X_i^{(j+1)}$ for all $i$ and $j$, we take 
$F^{(p)}:=D^p(F)$. The jet scheme $Z_m$ is defined by the ideal
$(F,F',\ldots,F^{(m)})$.

For every $m\geq 1$, there are canonical projections
$\phi_m:W_m\longrightarrow W_{m-1}$ induced by the truncation
homomorphisms $\CC[t]/(t^{m+1})\longrightarrow\CC[t]/(t^m)$.
By composing these projections we get morphisms $\rho_m:W_m
\longrightarrow W$.

If $W$ is a smooth, connected variety, then $W_m$ is smooth, connected,
and $\dim\,W_m=(m+1)\dim\,W$, for all $m$. It follows from the definition
that taking jet schemes commutes with open immersions. In particular,
if $W$ has pure dimension $n$, then $\rho_m^{-1}(W_{\rm reg})$
is smooth, of pure dimension $(m+1)n$.

For future reference, we record here two lemmas.
We denote by $[\cdot]$ the integral part function.

\begin{lemma}\label{fiber}{\rm (\cite{Mu1} 3.7)}
If $X$ is a smooth, connected variety of dimension $n$, $D\subset X$
is an effective divisor, and $x\in D$ is a point with ${\rm mult}_xD=q$,
then 
$$\dim\rho_m^{-1}(x)\leq mn-[m/q],$$
for every $m\in\NN$.
\end{lemma}

In fact, the only assertion we will need from Lemma~\ref{fiber}
is that $\dim\,\rho_m^{-1}(x)\leq mn-1$, if $m\geq q$, which follows easily
from the equations describing the jet schemes. 

\smallskip

If we have a family of schemes $\pi:{\mathcal W}\longrightarrow 
S$, we denote the fiber $\pi^{-1}(s)$ by ${\mathcal W}_s$.
The projection morphism $({\mathcal W}_s)_m\longrightarrow 
{\mathcal W}_s$ will be denoted by $\rho_m^{{\mathcal W}_m}$.

\begin{lemma}\label{semicont}{\rm (\cite{Mu2} 2.3)}
Let $\pi:{\mathcal W}\longrightarrow S$ be a family of schemes
and $\tau:S\longrightarrow{\mathcal W}$ a section of $\pi$.
For every $m\in\NN$, the function
$$f(s)=\dim(\rho_m^{{\mathcal W}_s})^{-1}(\tau(s))$$
is upper semicontinuous on the set of closed points of $S$.
\end{lemma}

We give now the proofs of our results.

\begin{proof}[Proof of Theorem~\ref{main1}]
If $\rho_m:Z_m\longrightarrow Z$ is the canonical projection,
then we have an isomorphism
\begin{equation}\label{isom}
\rho_m^{-1}(0)\simeq Z_{m-d}\times\AA^{n(d-1)},
\end{equation}
for every $m\geq d-1$ (we put $Z_{-1}=\{0\}$).
Indeed, for a $\CC$-algebra $A$, an $A$-valued point
of $\rho_m^{-1}(0)$ is a ring homomorphism
$$\phi:\CC[X_1,\ldots,X_n]/(F)\longrightarrow
A[t]/(t^{m+1}),$$
such that $\phi(X_i)\in(t)$ for all $i$.
Here $F$ is an equation defining $Z$. Therefore we can write
$\phi(X_i)=tf_i$, and $\phi$ is a homomorphism if and only if
the classes of $f_i$ in $A[t]/(t^{m+1-d})$ define an
$A$-valued point of $Z_{m-d}$. But $\phi$ is uniquely 
determined by the classes of $f_i$ in $A[t]/(t^m)$, so this proves 
the isomorphism in equation~(\ref{isom}).

Recall that we have $\dim\,{\rm Sing}\,(Z)=r$.
 An easy application
of Lemma~\ref{semicont} shows that for every
$x\in{\rm Sing}(Z)$, we have $\dim\,\rho_m^{-1}(x)
\leq\dim\,\rho_m^{-1}(0)$. 
We deduce that if $m\geq d-1$, then
$$\dim\,Z_m\leq\max\{(m+1)(n-1), r+\dim\,\rho_m^{-1}(0)\}$$
$$=\max\{(m+1)(n-1),\dim\,Z_{m-d}+n(d-1)+r\}.$$
A recursive application of this inequality shows that
for every $p\geq 1$, we have $\dim\,Z_{pd-1}\leq pd(n-1)$,
if $d\leq n-r$, and $\dim\,Z_{pd-1}\leq p(nd-n+r)$,
if $d\geq n-r$. This implies
$$\dim\,Z_{pd-1}/pd\leq\max\{n-1,(nd-n+r)/d\}.$$

By Theorem~\ref{ingred}, there is $p\geq 1$ such that
$$c(\AA^n,Z)=n-\frac{\dim\,Z_{pd-1}}{pd},$$
and we get $c(\AA^n, Z)\geq\min\{1,(n-r)/d\}$.
\end{proof}

\smallskip

\begin{proof}[Proof of Theorem~\ref{main2}]
We know that $d\geq n-r+1$ and $c(\AA^n,Z)=(n-r)/d$. By Theorem~\ref{ingred},
there is $k\geq 1$ such that
\begin{equation}\label{for_k}
\dim\,Z_{kd-1}=k(nd-n+r).
\end{equation}

We first show that if $k\geq 2$ and equation~(\ref{for_k}) holds
for $k$, then it holds also for $k-1$. 

Since $c(\AA^n, Z)=(n-r)/d$, it follows from Theorem~\ref{ingred}
that
\begin{equation}\label{ineq1}
\dim\,Z_{(k-1)d-1}\leq (k-1)(nd-n+r).
\end{equation}
The isomorphism~(\ref{isom}) in the proof of Theorem~\ref{main1} implies
\begin{equation}\label{ineq2}
\dim\,\rho_{kd-1}^{-1}(0)\leq (k-1)(nd-n+r)+nd-n.
\end{equation}
On the other hand, we have
\begin{equation}\label{ineq3}
\dim\,Z_{kd-1}\leq\max\{kd(n-1), \dim\,\rho_{kd-1}^{-1}(0)+r\}.
\end{equation}

Since $kd(n-1)<k(nd-n+r)$, we deduce from (\ref{for_k}), (\ref{ineq2})
and (\ref{ineq3}) that we have equality in (\ref{ineq2}), hence in
(\ref{ineq1}). Therefore equation~(\ref{for_k}) holds also for $k-1$.

The above argument shows that  
equation~(\ref{for_k}) holds for $k=1$, so that we have
$$\dim\,Z_{d-1}=dn-n+r>d(n-1).$$
 Equation~(\ref{isom}) in the proof
of Theorem~\ref{main1} gives $\rho_{d-1}^{-1}(0)\simeq
\AA^{n(d-1)}$.

Since $\dim\,\rho_{d-1}^{-1}(x)\leq \dim\,\rho_{d-1}^{-1}(0)$
for every $x\in Z$,
we deduce that there is a closed subset $W\subseteq\,{\rm Sing}\,(Z)$
with $\dim\,W=r$ such that $\dim\,\rho_{d-1}^{-1}(x)
=\dim\,\rho_{d-1}^{-1}(0)$ for all $x\in W$.

Fix $x\in W$. If ${\rm mult}_xZ\leq d-1$, then
Lemma~\ref{fiber} would give $\dim\,\rho_{d-1}^{-1}(x)\leq (d-1)n-1$,
a contradiction. Therefore we must have ${\rm mult}_xZ\geq d$.

By a linear change of coordinates we may assume that $x=(1,0,\ldots,0)$
and we write the equation $F$ of $Z$ as
$F=\sum_{i=0}^df_i(X_2,\ldots,X_n)X_1^{d-i}$, with $f_i$ homogeneous of
degree $i$ for all $i$. Since 
$$\frac{\partial^{\,p}F}{\partial X_1^p}(1,0,\ldots,0)=0,$$
for $p<d$, we deduce that $f_i=0$ for $i<d$, so that $Z=T_0\times\AA^1$,
where $T_0$ is the hypersurface defined by $f_d$ in $\AA^{n-1}$.

Since $\dim\,W=r$, by applying inductively the above argument
we get a homogeneous hypersurface $T=T_{r-1}\subset\AA^{n-r}$
such that $Z=T\times\AA^r$. Note that our assumption on the singular locus
of $Z$ implies that ${\rm Sing}\,(Z)=\{0\}$.

The converse is standard, and in fact we will prove a slightly stronger
 statement in the next proposition.
\end{proof}

The following proposition is well-known, but we include a proof
for the benefit of the reader.

\begin{proposition}\label{standard}
Let $T\subset\AA^{n-r}$ be a hypersurface defined by a homogeneous
polynomial of degree $d\geq n-r$ such that ${\rm Sing}\,(T)=\{0\}$.
If $Z=T\times\AA^r\subset\AA^n$ and $Y\subset\PP^{n-1}$ is the 
projectivization of $Z$, then $c(\PP^{n-1},Y)=(n-r)/d$.
\end{proposition}

\begin{proof}
An equivalent statement with that of the proposition
is that $c(\AA^n\setminus\{0\},Z\setminus\{0\})=(n-r)/d$.
Since $Z=T\times\AA^r$, we have $Z_m\simeq T_m\times\AA^{mr}$.
This implies that every irreducible component of $Z_m$
dominates $\AA^r$, so that $c(\AA^n\setminus\{0\},Z\setminus\{0\})=
c(\AA^n,Z)=c(\AA^{n-r},T)$.

 The fact that this number is $(n-r)/d$
is well-known. To see this using jet schemes we can use
the analogue of equation~(\ref{isom}) in the proof of 
Theorem~\ref{main1} together with the fact that
$$\dim\,T_m=\max\{(m+1)(n-r-1),\dim\,(\rho_m^T)^{-1}(0)\},$$
where $\rho_m^T$ is the projection corresponding to $T$.
By induction we get $\dim\,T_{kd-1}=k(n-r)(d-1)$, for every
$k\geq 1$, and Theorem~\ref{ingred}
gives $c(\AA^{n-r},T)=(n-r)/d$.
\end{proof}

\providecommand{\bysame}{\leavevmode \hbox \o3em
{\hrulefill}\thinspace}

\end{document}